\documentclass{article}

\PassOptionsToPackage{numbers, compress}{natbib}

 \usepackage[dblblindworkshop,final]{neurips_2025}

\workshoptitle{ScaleOPT: GPU-Accelerated and Scalable Optimization}

\usepackage[utf8]{inputenc} % allow utf-8 input
\usepackage[T1]{fontenc}    % use 8-bit T1 fonts
\usepackage{hyperref}       % hyperlinks
\usepackage{url}            % simple URL typesetting
\usepackage{booktabs}       % professional-quality tables
\usepackage{amsfonts}       % blackboard math symbols
\usepackage{nicefrac}       % compact symbols for 1/2, etc.
\usepackage{microtype}      % microtypography
\usepackage{xcolor}         % colors
\usepackage{amsmath}
\usepackage{algorithm}
\usepackage{algpseudocode}
\usepackage{verbatim}

\title{GPU-Accelerated Primal Heuristics for Mixed Integer Programming}

\author{%
  Akif Çördük \\
  Nvidia \\
  \texttt{acoerduek@nvidia.com} \\
  \And
  Piotr Sielski \\
  Nvidia \\
  \texttt{psielski@nvidia.com} \\
  \AND
  Alice Boucher \\
  Nvidia \\
  \texttt{yboucher@nvidia.com} \\
  \And
  Kumar Aatish \\
  Nvidia \\
  \texttt{kaatish@nvidia.com} \\
}

\begin{document}

\maketitle

\begin{abstract}
  We introduce a fusion of GPU accelerated primal heuristics for Mixed Integer Programming. Leveraging GPU acceleration enables exploration of larger search regions and faster iterations. A GPU-accelerated PDLP serves as an approximate LP solver, while a new probing cache facilitates rapid roundings and early infeasibility detection. Several state-of-the-art heuristics, including Feasibility Pump, Feasibility Jump, and Fix-and-Propagate, are further accelerated and enhanced. The combined approach of these GPU-driven algorithms yields significant improvements over existing methods, both in the number of feasible solutions and the quality of objectives by achieving 221 feasible solutions and 22\% objective gap in the MIPLIB2017 benchmark on a presolved dataset.
\end{abstract}

\section{Introduction}

Mixed Integer Linear Programming (MILP) problem is a mathematical optimization problem defined as a Linear Programming problem with an additional requirement that some of the variables must be integers. MILP is a generalization of NP-hard \cite{Karp} Integer Programming and is defined as follows: 
\begin{align*}
\text{Minimize} \quad & \sum_{i=1}^p c_i x_i + \sum_{j=p+1}^n d_j y_j \\
\text{subject to} \quad & \sum_{i=1}^p a_{ki} x_i + \sum_{j=p+1}^n b_{kj} y_j \leq g_k, \quad \forall k \in \{1,\dots,m\} \\
& x_i \in [l_i,u_i], \quad \forall i \in \{1,\dots,p\} \quad
y_j \in [l_j,u_j], \quad \forall j \in \{p+1,\dots,n\} \\
& y_j \in \mathbb{R}, \quad x_i \in \mathbb{Z} \quad
 c_i, d_j, a_{ki}, b_{kj}, g_k \in \mathbb{R}
\end{align*}

Currently, the state-of-the-art method for solving MILP problems to optimality is the Branch and Cut which works by enumerating the tree of all possibilities (Branching) and strengthening the formulation by adding cutting planes at tree nodes. During the process of tree exploration, many branches can be fathomed. One of the possibilities of fathoming is by comparing the objective of the best incumbent solution with the dual bound of the branch. Such a solution may be obtained by using various primal heuristic methods- methods that focus on finding strong feasible solutions without a proof of optimality. 

In many cases the MILP problem is too hard to solve to optimality in a given time budget, in which case the primal heuristics 
serve as one of the main tools for finding good feasible solutions. Significant efforts have been made in recent years to
improve such methods \cite{FixPropagate2025,FJ,FJ20,FPwObjective,FP20,GeneralFP,FP10years}. Feasibility Pump(FP) is the most commonly heuristics that uses sequence of projections and roundings. The projections minimize the L1 distance of integer variables to the polytope which has a converging behavior that the sum of L1 distances decrease over iterations. When the distance of the projection is zero, an integer feasible solution is found. Often, projection and rounding enter a cycle that is broken by various methods, including perturbation  \cite{FischettiGloverLodi2005} and the WalkSat algorithm  \cite{Dey2018,SelmanKautzCohen1994}. Various versions of FP has been proposed, some of which use objective components \cite{FPwObjective} and two-stage runs starting from binary variables first by relaxing the other integer variables \cite{GeneralFP}. As rounding heuristics, one of the state-of-the-art methods is the fix-and-propagate  \cite{FixPropagate2025} which does a sequence of variable fixing (rounding the variable and fixing its bounds to the rounding) and bounds propagation(BP). The authors handle the infeasibility in two ways, backtracking and repair procedure which does a neighborhood search by shifting the variable bounds to minimize constraint violation. 

We propose a fusion of methods that work together to achieve a higher number of feasible solutions and better objective than the state-of-the-art methods. The proposed framework, which is an improved FP, includes strategic use of the Local-MIP\cite{FJ20} to break cycles in the FP and also to search a neighborhood of an integer point. We also improve fix-and-propagate by introducing a new method of bulk rounding which uses a probing cache. In addition, we provide efficient GPU implementation of bound propagation algorithm with double probing, load balancing and propagation on changed constraints. We also introduce variable rounding order prioritization on new criterion. The PDLP algorithm is used in the beginning and in FP projections with approximate results and warm starts \cite{nicolasblog,Applegate2021,mexi2023scylla} that allows faster iterations. All methods are parallelized and implemented on Nvidia CUDA. The whole framework and the code base is open source and accessible at Nvidia cuOpt Github repository \cite{cuOpt}.

\section{GPU Algorithms}

\subsection{Bound Propagation} 

The bound propagation(BP) \cite{BPSavelsbergh} is a two step iterative process for tightening the bounds of variables. In every iteration, the minimal and maximal constraint activities are calculated for every constraint $k$ :

\begin{align*}
\sum_{i=1}^n a_{ki} z_i \leq g_k \quad \text{where} \quad
& z_{i} \in [l_{i}, u_{i}] \quad
L^{k}_{min} = \sum_{\substack{i \in \mathbb{N} \\ a_{ki} > 0}} a_{ki} l_i + \sum_{\substack{i \in \mathbb{N} \\ a_{ki} < 0}} a_{ki} u_i
\end{align*}

The activity values, along with the current variable bounds are used to tighten the bounds further for every variable $i$ :

$$
\hat{l_i} = \max_{k \in \mathbb{N}}\{\frac{g_k - L^{k}_{min} + a_{ki}u_i}{a_{ki}} : a_{ki} < 0\}, \quad \hat{u_i} = \min_{k \in \mathbb{N}}\{\frac{g_k - L^{k}_{min} + a_{ki}l_i}{a_{ki}} : a_{ki} > 0\}
$$

If the bound updates do not cross a threshold, the iterations are terminated.

In the MIPLIB2017 \cite{MIPLIB2017} instances, the median of average number of terms per constraint or variable was found to be in the range 6-9. We also see a substantial variability in the number of terms across the dataset. This observation motivates us to specialize our approach tailored to the number of non-zero terms for a constraint or a variable similar to \cite{SOFRANAC2022102874} but we go a few steps further.

We create a load balanced version of the problem statement using Logarithmic Radix Binning\cite{LRB}. Constraints with similar number of non-zero terms are binned together in global memory to form a Compressed Sparse Row(CSR). A similar operation is done for each variable to obtain a load balanced version of the CSR transpose.

The activity and bounds update calculation is done via three different approaches.
For activity calculation, the sub-warp approach splits each warp to work on multiple small constraints at once. Each warp accesses the binning offsets generated by LRB and determines the bin assigned to it. This information, along with the warp id helps in calculating which constraints the warp is operating on and how the warp should be split to handle the multiple constraints at once. This approach allows us to work on all the coefficient and variable index data brought in to the cache.
The block approach works on medium sized constraints where a single thread block is sufficient to operate on a constraint.
The third, heavy element approach splits every large constraint (greater than 16k terms) to be dealt by multiple blocks. The partial summations of activities calculated by each block are written to a global temporary array. A second pass summation is then run to reach the final activity values of the heavy constraints.
A similar approach is utilized for bounds update step, where the summation is replaced with minimum and maximum operators.

In both the steps while the CSR information is loaded into memory in a coalesced fashion, random accesses to two values for each non-zero term is necessary, i.e. in the activity calculation requires variable lower and upper bound, and bounds update requires the minimal and maximal activity. In order to improve cache performance, the variable bounds are grouped together into a single buffer so that the lower and upper bounds of variable $i$ are at location $2i$ and $2i+1$. A single buffer is maintained for the activities as well.

CUDA graphs have been utilized to reduce kernel launch latency and run the three approach concurrently.

\subsection{Probing Cache} \label{probecache}

For a binary variable, the implications derived from fixing its value to $0$ or $1$ can be computed beforehand and reused later in other procedures. For the general integer variables, computing implications for all possible values is not practically possible, so we only search for implication of partitioning to two sub-intervals: For boxed variables, bound implications for $x < (u + l)/2$ and $x\geq (u + l)/2$ are precomputed. For one-sided bounds, say bounded by lower bound $l$, we consider the following possibilities: $x = l$ and $x \geq l+1$ . We store the implications in a cache which is used to warm-start bound propagation and create promising insights on variable value combinations. The key place where probing cache is used is the rounding heuristics, where we prepare a bulk by analyzing the conflicts and the resulting cached bounds are used as a warm start for the bound propagation. This allows large speedups in BP rounding, and as a result we can run the procedure more often. 

The computation of the probing cache might be too expensive to compute for all variables, so prioritization is needed. We prioritize variables according to 3 criteria lexicographically: total number of violated constraints per variable, maximum violation caused by a variable, minimum unit slack across constraints. The first two criteria are enforced on variables that appear in different constraints with negated signs. The minimum activity contribution in one constraint will not be valid in another constraint (the activity will consider lower bound in one constraint and upper bound in another constraint), meaning that these variables are more likely to cause conflicts during probing. On these variables, we compute the activity on each constraints by flipping the variable contribution from lower bound to upper bound and visa versa. This reveals possible constraint violations when a variable is forced to take extreme values of its bounds. For each variable, with this criteria, we compute the total number of violated constraints and maximum violation. The last sorting criteria is the unit slack consumption which is explained in \ref{rounding}.

The GPU implementation of BP might not utilize the GPU fully on single-variable propagations, as the number of constraints that include the variable might be small. For that reason, we launch multiple propagation procedures on different threads and cuda streams concurrently. This increases the utilization of the GPU by concurrent running of kernels and memory transfers and increases the number of variables cached in a given time limit. We also improve the performance by doing a double propagation for a variable: we propagate both ranges of the cache together with marginal overhead as majority of the memory accesses and computations are common. 

\subsection{Rounding Heuristics} \label{rounding}

In a heuristic algorithm, it is often needed to round a fractional solution. Some rounding methods consider variable bound propagation \cite{FP20,FixPropagate2025} which was improved and accelerated in this work. The procedure fixes a variable to an integer value and propagates its bounds to other variables at each iteration. The rounding candidate is chosen according to a distribution \cite{GeneralFP} which does random rounding for fractions closer to .5 and does a nearest rounding for fractions away from .5. Given that the number of variables can be large, the BP procedure could take a lot of time. The procedure is improved by bulk rounding that rounds multiple variables together quickly. We use the probing cache, generated at the beginning of the whole procedure (see \ref{probecache}), to decide which values are not conflicting and should be rounded together in the bulk. BP procedure is only run after a bulk is fixed. If it fails, we backtrack and run the BP for each variable for this bulk separately for each variable. During backtracking if we encounter an infeasible propagation, we run the repair procedure \cite{FixPropagate2025} to do a search for feasible bounds given fixed variables. The procedure tries to remove infeasibility from a constraint  by shifting the variable bounds one at a time. We compute a shift value for each variable in the constraint which makes the constraint feasible. The implementation of the repair procedure is parallelized on GPU. Even if the repair procedure makes the bounds feasible, the bounds could still be eventually infeasible after running the propagation procedure. Additionally, the modified bounds of the unfixed variables might be too tight or too relaxed as they are not generated by the propagation from fixed variables which is the ground truth for bounds \cite{FixPropagate2025}. The repair and propagation on original bounds sequence is run until the projection yields feasible bounds or we exhaust iteration or time limit. If the solution is deemed infeasible after we tried certain number of repair procedures, we continue bulk rounding until the end without any backtracking and repair to reduce the infeasibility of the final solution.

At each bulk propagation step, we generate two bulks which we propagate in parallel at marginal cost thanks to data locality. Both of the bulks are generated by drawing a value for a variable two times from the rounding distribution. For the fractional values that are closer to .5, we will likely have two different values and less fractional values will be represented as nearest rounding in the bulk. This protects the original rounding candidates while introducing some randomization for variables that could be rounded to both sides. The bulk size is kept proportional to the square root of the number of remaining variables. The last 36 variables are fixed and propagated one by one. After achieving a feasible propagation for all variables, we run a polishing LP to find feasible values for the continuous variables.

The order of propagation plays a crucial role in which value a variable will be fixed. We do an initial sort at the beginning and an implied slack sort (stable sort) at each bulk rounding iteration. The initial sort acts as a tie-breaker for the stable sort before each bulk rounding. In the initial sort, we consider variable bounds intervals in which binary and ternary variables always precede the other variables. Then, within binary, ternary and rest of the variables we sort according to the fractions of the variables. The implied slack sort dynamically sorts the remaining variables at each bulk. For each variable, we quantify how much a unit change in its value reduces the constraint slack. The constraint activity slack is calculated as $cnstr\_upper - min\_act$ or $cnstr\_lower - max\_act$. The variable unit slack consumption is calculated as $s = (a / act)$ where $a$ is the coefficient, $act$ is the relevant activity slack value. We calculate the total impact of slack consumptions across each constraint and then sort according to:
$
S_i = \sum_{c \in C} s_{ic}^2
$
where $S$ is the final sort value for var $i$, $C$ is the set of constraints that include the variable. A pseudo code is given in appendix \ref{alg:rounding}. 

\subsection{Local Search}
Recent advances in Mixed Integer Programming have capitalized on local search techniques to efficiently discover feasible and high-quality solutions, especially in the case of large instances or time-constrained scenarios.

A prominent example is the Feasibility Jump (FJ) heuristic \cite{FJ} which won the MIP 2022 Computational Competition, and is a primal heuristic designed to quickly find feasible solutions by iteratively minimizing constraint violations using a Lagrangian relaxation approach. However, it mainly focuses on feasibility, and lacks efficient mechanisms for objective chasing. Local-MIP \cite{FJ20} addresses these gaps in several ways. Two new specialized move operators are introduced:

The breakthrough move operator is designed to seek aggressive improvements in the objective value once feasibility is achieved, by targeting promising variable changes that might escape the current local minima.
The lift move operator exploits slack in satisfied constraints to improve the objective while maintaining feasibility. 
Local-MIP also adopts a hierarchical scoring function: at the first level, progress scores reward immediate feasibility and objective improvements, while the second level accounts for substantial objective change or satisfying a new constraint. 
Experimental results show that Local-MIP finds feasible solutions more than FJ and with better objective values. 

We have implemented a GPU-accelerated heuristic that builds directly upon these papers while introducing key improvements for parallel execution.

GPU-Local-MIP maintains the same principles of Lagrangian-based relaxation of linear constraints and leveraging violation penalty weights. However, it differs in implementation strategy and algo-technical choices:
Multiple candidate moves and their full scores are computed once per iteration. This is done in parallel, to saturate the GPU by processing a large number of constraint/variable pairs concurrently.
Move generation, scoring, and selection are decoupled into pipelined kernel stages via CUDA streams.

As in the original paper, GPU-Local-MIP supports three base move types: Mixed Tight Moves (MTM) that satisfy constraints by targeting bound violations, lift moves for objective improvement while maintaining feasibility, and breakthrough moves to pursue objective improvements in newly feasible regions.
The original FJ updated only the selected variable's score and approximated others using deltas, a local strategy that could inhibit global progress. GPU-Local-MIP instead recalculates all affected variable scores when needed on each iteration. For this purpose, an array of “neighbor” variables (corresponding to the distance-2 neighborhood of the problem bipartite graph, i.e. the variables that share a constraint with the selected variable) is precomputed during problem initialization using a sparse SpMM operation. This is practical for small to medium instances, but becomes unsuitable for large instances because the number of edges of a distance-2 graph is $O(N^2)$ with $N$ the number of nodes. To fit in available memory, only the neighborhood of the selected variable is maintained, and computed on the fly during each iteration.
Since MTM variable moves outside this set will not be affected by any resulting change in constraint LHS values from applying the selected move, they do not need to be recomputed.
Additionally, binary variables receive streamlined flip-based computation through an optimized dedicated code path and kernels.

A naive implementation mapping thread blocks to problem rows would be vulnerable to load-imbalance, with blocks mapped to smaller rows completing sooner than longer rows and leaving GPU resources underutilized. Rather than assigning one thread block per row or one warp per variable, GPU-Local-MIP uses a prefix-sum-based \cite{lbpsum} load balancing scheme that maps warps to groups of 32 non-zeros. This is maintained as a separate code-path, which is dynamically selected according to the size and the level of imbalance of the problem matrix. 

To improve numerical stability, Local-MIP integrates Kahan summation to update LHS values and constraint violations at each iteration. These are recomputed after a fixed number of iterations with a SpMV operation to minimize running error.

Candidate moves are filtered via a stream compaction step and passed into reduction kernels to select the highest scoring move.
GPU-Local-MIP utilizes CUDA Graphs to eliminate host kernel launch overhead.
Additionally, CUDA cooperative groups are used to ensure full-grid synchronizations where needed. 

\subsection{Extended Feasibility Pump}

We use a variety of heuristics in collaboration to powerfully utilize their strengths. The algorithm starts with LS run on a solution where we round the LP relaxation solution, which might be approximate, with the rounding heuristics. This is to find any easy feasible solution around the LP-optimal region. If no feasible is found, we run Local-MIP on an all-zero solution to exploit problem structures where majority of the variables are zero. 

The algorithms are fused around the FP algorithm \cite{FischettiGloverLodi2005} which includes projection of integer variables onto the polytope and rounding of this projection. Various versions of the FP exist \cite{GeneralFP,FPwObjective} including two-stage FP and objective FP which we utilize. In the FP, cycles of integer solutions or projection distance occur and various approaches have been proposed to break the cycles \cite{FP10years} mostly using random perturbations. The approach of using WalkSAT \cite{Dey2018} is a notable exception, but lacks generalization to integer variables. Also this approach only reduces the number of infeasible constraints, but does not take into account the total violation or hardness of the constraints. 
Instead we propose to break the cycles with Local-MIP which considers the general problem structure with lagrangian weights and is a strong heuristic by itself.

We start the fused heuristic algorithm on the LS solution which has a minimized weighted constraint violation. A sequence of FP projection, rounding heuristics and LS are run until the time limit is reached. The FP projection brings the search region closer to the polytope and helps bring continuous variables to the feasible region. The projection time is limited to 1 second for all problems. This allows us to obtain LP-optimal solutions for smaller/easier problems and approximate solutions, which can be feasible or infeasible, for larger/harder problems. Approximate solutions are still useful \cite{Applegate2021} as they are close to the polytope and provide promising rounding points. We warm start the PDLP algorithm with the previous primal and dual solutions of the previous projection. The added variables from FP projection can be warm started with the L1 distances from previous projection however the dual variables for the added constraints are not warm started. The warm start property works cumulatively over FP iterations and help make the solution feasible over time. This way, we can run more iterations without sacrificing too much feasibility.

At each FP projection, majority of the variables keep their integer values. Rounding heuristics are used to round the remaining fractional variables to an integer value. To avoid inevitable infeasibility resulting from the starting integer values, we round all variables from the beginning. This allows three things: utilizing variable order coming from sorting, avoiding infeasibility by moving the already integer value to a new feasible value and utilizing the repair procedure to move the value dynamically across different stages of rounding. The Local-MIP algorithm is used after each rounding to test whether there is a feasible solution around the rounded point by running it with 20\% of the and rounding time. If a feasible solution is found we break the algorithm, otherwise we discard the solution to continue FP trajectory uninterrupted. The rounding heuristics disturb the converging property of FP which always guarantees lower or equal L1 distance. To detect cycles, we use a method similar to in \cite{FPwObjective} where we detect the cycle when projection distance is not improved by more than 10\% of the average of last $N$ projections.

The whole procedure can be extended to be used as an improvement heuristic rather than just generation heuristic. At each feasible point, we add a cutting plane perpendicular to the objective direction as upper bound \cite{FischettiGloverLodi2005}. Thus, the algorithm works on a restricted problem in which any new feasible point would be better than our previous solution. A pseudo code can be found in appendix \ref{alg:fused_heuristic}.

\section{Results}

We compare our heuristics, which was run with the $fp\_only=true$ flag in Nvidia cuOpt \cite{cuOpt}, with the results obtained in the Local-MIP\cite{FJ20} paper and the Fix-And-Propagate \cite{FixPropagate2025} paper.
We employ a similar protocol to \cite{FixPropagate2025}: each method is run with a time limit of 10 minutes on the 240 instances of the MIPLIB2017 \cite{MIPLIB2017} benchmark set, using three separate seeds and averaging out the results. Each instance has been presolved using the open-source solver Highs \cite{highs} before each run in order to match the results reported by \cite{FixPropagate2025} which used the commercial solver CPLEX for this purpose. Our GPU method was executed on a H100 80GB HBM3 GPU coupled with a Intel(R) Xeon(R) Platinum 8480CL E5 CPU.
The Local-MIP method was run on a AMD EPYC 7742 CPU \footnote{Using the open-source implementation provided at \url{https://github.com/linpeng0105/Local-MIP}}, commit hash 141a2f4, with the addition of code to support "RANGES"-type constraints. 
The feasibility tolerance used is an absolute tolerance of 1e-6 on all experiments and the references.

In Table \ref{tab:results-miplib2017}, we present the results of a 10 min run in terms of number of feasible solutions and primal gap w.r.t. to the optimal solution, which we define as such:
\[
\text{Primal gap} =
\begin{cases}
1, & \text{if } z_{\text{obj}} \cdot z_{\text{optimal}} < 0 \text{ or infeasible}, \\[6pt]
\frac{|z_{\text{obj}} - z_{\text{optimal}}|}{\max\left\{ |z_{\text{obj}}|, |z_{\text{optimal}}| \right\}}, & \text{otherwise}.
\end{cases}
\]

\begin{table}[ht]
\centering
\begin{tabular}{lccc}
% \toprule
Method    & Avg. \#Feasible & Primal gap \\
\midrule
Fix-And-Propagate portfolio default & 193.8     & 0.66    \\
Local-MIP & 188.67 & 0.46 \\
GPU Local-MIP & 205 & 0.41 \\
GPU Extended FP with Nearest Rounding & 220 & 0.23 \\
GPU Extended FP with Fix and Propagate & 220.67 & 0.22 \\
\bottomrule
\end{tabular}
\caption{Feasibility \& primal gap results over the presolved MIPLIB2017 set.}
\label{tab:results-miplib2017}
\end{table}

GPU-accelerated Local-MIP outperforms the CPU version by both number of feasible solutions and average primal gap. The extended FP framework outperforms the state-of-the-art references significantly even with the base version of FP with the nearest rounding. The repair procedure does not contribute strongly to the results in our tests, so it was disabled during the experiments. The fix-and-propagate rounding adds slightly more feasible solutions 0.67 and reduces 1\% objective.

\medskip

\bibliographystyle{abbrvnat}
\bibliography{main}

%%%%%%%%%%%%%%%%%%%%%%%%%%%%%%%%%%%%%%%%%%%%%%%%%%%%%%%%%%%%

\appendix

\section{Technical Appendices and Supplementary Material}

\begin{algorithm}[H]
\caption{Bound Propagation Rounding}
\label{alg:rounding}
\begin{algorithmic}
\Procedure{PROPAGATION-ROUND}{$sol, timer$}
    \State Initialize \textit{unset\_vars} with all integer variables from $sol$.
    \State Initialize control flags: \textit{rounding\_infeasible}, \textit{recovery\_mode} $\gets \text{false}$.

    \State \Call{Sort-By-Interval-And-Fractionality}{$sol, \text{\textit{unset\_vars}}$}.

    \While{\textit{set\_count} $<$ \textit{unset\_vars}.size() and not \Call{Time-Limit-Reached}{}}
        \State \textit{bulk\_size} $\gets$ \Call{Get-Bulk-Size}{\textit{unset\_vars}.size()}.
        \Comment{Returns 1 if in recovery\_mode}

        \If{\textit{bulk\_size} $>$ 1}
            \State \Call{Sort-By-Implied-Slack}{\text{\textit{unset\_vars}}} \Comment{Stable sort}
        \EndIf

        \State \textit{vars\_to\_set} $\gets$ next \textit{bulk\_size} variables from \textit{unset\_vars}.
        \State \textit{probe\_vec\_0}, \textit{probe\_vec\_1} $\gets$ \Call{Generate-Candidate-Values}{$sol, \text{\textit{vars\_to\_set}}$}.
        \State \textit{results} $\gets$ \Call{Parallel-Propagate}{$sol, \text{\textit{probe\_vec\_0}}, \text{\textit{probe\_vec\_1}}$} 
        \State \textit{infeas\_count\_0}, \textit{infeas\_count\_1} $\gets$ \textit{results.infeas\_counts}.

        \If{\textit{infeas\_count\_0} == 0 or \textit{infeas\_count\_1} == 0} \Comment{At least one probe was successful.}

            \State \Call{Update-Solution-Bounds-With-Feasible}{$sol, \text{\textit{results}}, \text{\textit{selected\_update}}$}.
            \State \textit{set\_count} $\gets$ \Call{Count-All-Fixed-Vars}{$sol, \text{\textit{unset\_vars}}$}.
            \State \textit{recovery\_mode} $\gets \text{false}$.
        \Else \Comment{Both probes failed; trigger recovery or repair.}
            \If{not \textit{recovery\_mode} and not \textit{rounding\_infeasible}}
                \State \textit{recovery\_mode} $\gets \text{true}$. \Comment{Backtrack and switch to single-variable mode.}
            \Else
                \State \textit{rounding\_infeasible} $\gets \text{true}$. \Comment{Mark path as infeasible.}
            \EndIf
        \EndIf

        \If{\textit{rounding\_infeasible}}
            \State \textit{repaired} $\gets$ \Call{Run-Repair-Procedure}{$sol,\text{timer}$}.
            \If{\textit{repaired}}
                \State \textit{rounding\_infeasible} $\gets \text{false}$. \Comment{Reset flags after successful repair.}
            \EndIf
        \EndIf
    \EndWhile

    \If{final solution is feasible}
        \State \Call{Run-LP-Polish}{$sol$}.
    \EndIf
    \State \Return $sol$
\EndProcedure
\end{algorithmic}
\end{algorithm}

\begin{algorithm}
\caption{Fused Heuristic for Mixed-Integer Programming}
\label{alg:fused_heuristic}
\begin{algorithmic}[1]
\Procedure{FusedHeuristic}{$P, T_{limit}$}
   
    \State $s_{LP} \gets \text{SolveLPRelaxation}(P)$  \Comment{Initialize with a solution near the LP-optimal region}
    \State $s \gets \text{RoundingHeuristic}(s_{LP})$
    \State $(s_{best}, \text{feasible}) \gets \text{LocalSearch}(s)$
    
    \If{\textbf{not} feasible}
        \State $s_0 \gets \text{AllZeroSolution}()$ \Comment{If no solution found, try Local-MIP on an all-zero solution}
        \State $(s_{best}, \text{feasible}) \gets \text{EvolutionaryLocalSearch}(s_0)$
    \EndIf
    
    \State $s \gets s_{best}$ \Comment{Start main loop with the best solution found so far}
    \State $proj\_history \gets \emptyset$
    
    \While{time elapsed $< T_{limit}$}
        \If{feasible}
            \State \Comment{Improvement Heuristic: Add objective cutting plane}
            \State $P \gets P \cup \{\text{objective cut based on } s_{best}\}$ 
        \EndIf
        
        \State \Comment{Project integer solution onto the polytope (Feasibility Pump)}
        \State $s_{proj} \gets \text{FP\_Project}(s, P, \text{1s\_limit})$ \Comment{Warm-start with previous solutions}
        \State Add $L1\_distance(s, s_{proj})$ to $proj\_history$
        
        \State 
        \State $s_{rounded} \gets \text{RoundAllVars}(s_{proj})$ \Comment{Round all variables from the projected point}
        
        \State \Comment{Check for a feasible solution around the rounded point}
        \State $(s_{test}, \text{is\_new\_sol}) \gets \text{Local-MIP}(s_{rounded}, 0.2 \times T_{rounding})$
        
        \If{is\_new\_sol}
            \State $s_{best} \gets s_{test}$
            \State \textbf{break} \Comment{New feasible solution found, exit}
        \EndIf
        
        \State \Comment{Detect and break cycles in the FP trajectory}
        \If{CycleDetected($proj\_history$)}
            \State \Comment{Use Local-MIP to escape the cycle by minimizing weighted violation}
            \State $s \gets \text{Local-MIP}(s, \text{lagrangian\_weights})$
        \Else
            \State \Comment{Continue FP trajectory from the last integer solution}
            \State $s \gets \text{RoundToInteger}(s_{proj})$
        \EndIf
    \EndWhile
    
    \State \Return $s_{best}$
\EndProcedure
\end{algorithmic}
\end{algorithm}

\end{document}